\documentclass[twoside,10pt]{article}
\usepackage{amsmath,amssymb,amsthm}

\begin{document}

\setlength{\textwidth}{126mm} \setlength{\textheight}{180mm}
\setlength{\parindent}{0mm} \setlength{\parskip}{2pt plus 2pt}

\frenchspacing

\pagestyle{myheadings}

\markboth{Dimitar Mekerov}{Lie groups as $4$-dimensional
nonintegrable almost product manifolds}

\newtheorem{thm}{Theorem}[section]
\newtheorem{lem}[thm]{Lemma}
\newtheorem{prop}[thm]{Proposition}
\newtheorem{cor}[thm]{Corollary}
\newtheorem{probl}[thm]{Problem}

\newtheorem{defn}{Definition}[section]
\newtheorem{rem}{Remark}[section]
\newtheorem{exa}{Example}



\newcommand{\X}{\mathfrak{X}}
\newcommand{\B}{\mathcal{B}}
\newcommand{\s}{\mathfrak{S}}
\newcommand{\g}{\mathfrak{g}}
\newcommand{\W}{\mathcal{W}}
\newcommand{\Lgr}{\mathrm{L}}
\newcommand{\dd}{\mathrm{d}}

\newcommand{\pd}{\partial}
\newcommand{\ddx}{\frac{\pd}{\pd x^i}}
\newcommand{\ddy}{\frac{\pd}{\pd y^i}}
\newcommand{\ddu}{\frac{\pd}{\pd u^i}}
\newcommand{\ddv}{\frac{\pd}{\pd v^i}}

\newcommand{\diag}{\mathrm{diag}}
\newcommand{\End}{\mathrm{End}}
\newcommand{\im}{\mathrm{Im}}
\newcommand{\id}{\mathrm{id}}
\newcommand{\ad}{\mathrm{ad}}

\newcommand{\ie}{i.e.}
\newfont{\w}{msbm9 scaled\magstep1}
\def\R{\mbox{\w R}}
\newcommand{\norm}[1]{\left\Vert#1\right\Vert ^2}
\newcommand{\nN}{\norm{N}}
\newcommand{\nP}{\norm{\nabla P}}
\newcommand{\tr}{{\rm tr}}

\newcommand{\nJ}[1]{\norm{\nabla J_{#1}}}
\newcommand{\thmref}[1]{Theorem~\ref{#1}}
\newcommand{\propref}[1]{Proposition~\ref{#1}}
\newcommand{\secref}[1]{\S\ref{#1}}
\newcommand{\lemref}[1]{Lemma~\ref{#1}}
\newcommand{\dfnref}[1]{Definition~\ref{#1}}

\frenchspacing


\title{Lie groups as $4$-dimensional Riemannian
or pseudo-Riemannian almost product manifolds
with nonintegrable structure}

\author{Dimitar Mekerov}

\maketitle

{\small
{\it Abstract.} A Lie group as a 4-dimensional pseudo-Riemannian
manifold is considered. This manifold is equipped with an almost
product structure and a Killing metric in two ways. In the first
case à Riemannian almost product manifold with nonintegrable
structure is obtained, and in the second case -- a
pseudo-Riemannian one. Each belongs to a 4-parametric family of
manifolds, which are characterized geometrically.

}



\section{Preliminaries}

Let $M$ be a differentiable manifold with a tensor field $P$ of
type $(1,1)$ and a Riemannian metric $g$ such that
\begin{equation}\label{1.1}
    P^2=id,\quad g(Px,Py)=g(x,y)
\end{equation}
for arbitrary $x$, $y$ of the algebra $\X(M)$ of the smooth vector
fields on $M$. The tensor field $P$ is called an \emph{almost
product structure}. The manifold $(M,P,g)$ is called a
\emph{Riemannian} (\emph{pseudo-Riemannian}, resp.) \emph{almost
product manifold}, if $g$ is a Riemannian (pseudo-Riemannian,
resp.) metric. If $\tr{P}=0$, then $(M,P,g)$ is an
even-dimensional manifold. The classification from
\cite{StGr:connect} of Riemannian almost product manifolds is made
with respect to the tensor field $F$ of type (0,3), defined by
\begin{equation}\label{1.2}
F(x,y,z)=g\left(\left(\nabla_x P\right)y,z\right),
\end{equation}
where $\nabla$ is the Levi-Civita connection of $g$. The tensor
$F$ has the following properties:
\[
    F(x,y,z)=F(x,z,y)=-F(x,Py,Pz),\quad F(x,y,Pz)=-F(x,Py,z).
\]
In the case when $g$ is a pseudo-Riemannian metric, the same
classification is valid for pseudo-Riemannian almost product
manifolds, too. In these classifications the condition
\begin{equation}\label{sigma}
    F(x,y,z)+F(y,z,x)+F(z,x,y)=0
\end{equation}
defines a class $\W_3$, which is only the class of the three basic
classes $\W_1$, $\W_2$ and $\W_3$ with nonintegrable structure
$P$.

The class $\W_0$, defined by the condition $F(x,y,z)=0$, is
contained in the other classes. For this class $\nabla P=0$ and
therefore it is an analogue of the class of K\"ahlerian manifolds
in the almost Hermitian geometry.

The curvature tensor field $R$ is defined by $R(x,y)z=\nabla_x
\nabla_y z - \nabla_y \nabla_x z - \nabla_{[x,y]}z$ and the
corresponding tensor field of type $(0,4)$ is determined by
$R(x,y,z,w)=g(R(x,y)z,w)$.

Let $\{e_i\}$ be a basis of the tangent space $T_pM$ at a point
$p\in M$ and $g^{ij}$ be the components of the inverse matrix of
$g$ with respect to $\{e_i\}$. Then the Ricci tensor $\rho$ and
the scalar curvature $\tau$ are defined as follows
\begin{equation}\label{1.3}
    \rho(y,z)=g^{ij}R(e_i,y,z,e_j),
\end{equation}
\begin{equation}\label{1.4}
    \tau=g^{ij}\rho(e_i,e_j).
\end{equation}
The square norm of $\nabla P$ is defined by
\begin{equation}\label{1.5}
    \nP=g^{ij}g^{ks}g\left(\left(\nabla_{e_i}P\right)e_k,\left(\nabla_{e_j}P\right)e_s\right).
\end{equation}
It is clear that $\nabla P=0$ implies $\nP=0$ but the inverse
implication for the pseudo-Riemannian case is not always true. We
shall call a pseudo-Riemannian almost product manifold
\emph{isotropic $P$-manifold} if $\nabla P=0$.

The Weyl tensor on a $2n$-dimensional pseudo-Riemannian manifold
($n\geq 2$) is
\begin{equation}\label{1.6}
    W=R-\frac{1}{2n-2}\left(\psi_1(\rho)-\frac{\tau}{2n-1}\pi_1\right),
\end{equation}
where
\[
    \begin{array}{l}
      \psi_1(\rho)(x,y,z,w)=g(y,z)\rho(x,w)-g(x,z)\rho(y,w) \\
      \phantom{\psi_1(\rho)(x,y,z,w)}+\rho(y,z)g(x,w)-\rho(x,z)g(y,w); \\
      \pi_1(x,y,z,w)=g(y,z)g(x,w)-g(x,z)g(y,w). \\
    \end{array}
\]

Moreover, for $n\geq 2$ the Weyl tensor $W$ is zero if and only if
the manifold is \emph{conformally flat}.

If $\alpha$ is a non-degenerate 2-plane spanned by vectors $x, y
\in T_pM, p\in M$, then its sectional curvature is
\begin{equation}\label{1.7}
    k(\alpha)=\frac{R(x,y,y,x)}{\pi_1(x,y,y,x)}.
\end{equation}


\section{A Lie group as a 4-dimensional
pseudo-Rie\-mannian manifold with Killing metric}

Let $V$ be a real 4-dimensional vector space with a basis
$\{E_i\}$. Let us consider a structure of a Lie algebra determined
by commutators $[E_i,E_j]=C_{ij}^k E_k$, where $C_{ij}^k$ are
structure constants satisfying the anti-commutativity condition
$C_{ij}^k=-C_{ji}^k$ and the Jacobi identity $C_{ij}^k
C_{ks}^l+C_{js}^k C_{ki}^l+C_{si}^k C_{kj}^l=0$.

Let $G$ be the associated connected Lie group and $\{X_i\}$ be a
global basis of left invariant vector fields which is induced by
the basis $\{E_i\}$ of $V$. Then we have the decomposition
\begin{equation}\label{2.1}
    [X_i,X_j]=C_{ij}^k X_k.
\end{equation}

Let us consider the manifold $(G,g)$, where $g$ is a metric
determined by the conditions
\begin{equation}\label{2.2}
\begin{array}{c}
  g(X_1,X_1)=g(X_2,X_2)=g(X_3,X_3)=g(X_4,X_4)=1, \\[4pt]
  g(X_i,X_j)=0\quad \text{for}\quad i\neq j \\
\end{array}
\end{equation}
or by the conditions
\begin{equation}\label{2.3}
\begin{array}{c}
  g(X_1,X_1)=g(X_2,X_2)=-g(X_3,X_3)=-g(X_4,X_4)=1, \\[4pt]
  g(X_i,X_j)=0\quad \text{for}\quad i\neq j. \\
\end{array}
\end{equation}

Obviously, $g$ is a Riemannian metric if it is determined by
\eqref{2.2} and $g$ is a pseudo-Riemannian metric of signature
(2,2) if it is determined by \eqref{2.3}.

It is known that the metric $g$ on the group $G$ is called a
\emph{Killing metric} \cite{Hel} if the following condition is
valid
\begin{equation}\label{2.5}
    g([X,Y],Z)+g([X,Z],Y)=0.
\end{equation}
where $X$, $Y$, $Z$ are arbitrary vector fields.

If $g$ is a Killing metric, then according to the proof of
Theorem~2.1 in \cite{MaGrMe-4} the manifold  $(G,g)$ is
\emph{locally symmetric}, i.e. $\nabla R=0$. Moreover, the
components of $\nabla$ and $R$ are respectively
\begin{equation}\label{2.6}
    \nabla_{ij}=\nabla_{X_i} X_j=\frac{1}{2}[X_i,X_j],
\end{equation}
\begin{equation}\label{2.7}
    R_{ijks}=R(X_i,X_j,X_k,X_s)=-\frac{1}{4}g\left([X_i,X_j],[X_k,X_s]\right).
\end{equation}


\section{A Lie group as a Riemannian almost product manifold
with Killing metric and nonintegrable structure}

In this section we consider a Riemannian manifold $(G,P,g)$ with a
metric $g$ determined by \eqref{2.2} and a structure $P$ defined
as follows
\begin{equation}\label{3.1}
    PX_1=X_3,\quad PX_2=X_4,\quad PX_3=X_1,\quad PX_4=X_2.
\end{equation}
Obviously, $P^2=\id$. Moreover, \eqref{2.2} and \eqref{3.1} imply
\begin{equation}\label{3.2}
    g(PX_i,PX_j)=g(X_i,X_j).
\end{equation}
Therefore, $(G,P,g)$ is a Riemannian almost product manifold.

For the  manifold $(G,P,g)$ we propose that $g$ be a Killing
metric. Then $(G,P,g)$ is locally symmetric.

From \eqref{2.6} we obtain
\begin{equation}\label{3.3}
    \left( \nabla_{X_i} P \right)X_j=\frac{1}{2}\bigl([X_i,PX_j]-P[X_i,X_j]\bigr).
\end{equation}
Then, according to \eqref{1.2}, for the components of $F$ we have
\begin{equation}\label{3.4}
    F_{ijk}=\frac{1}{2}g\bigl([X_i,PX_j]-P[X_i,X_j],X_k\bigr).
\end{equation}
Hence, having in mind \eqref{2.2}, \eqref{3.1} and \eqref{3.2}, we
get
\begin{equation}\label{3.5}
    F_{ijk}+F_{jki}+F_{kij}=0,
\end{equation}
i.e. $(G,P,g)$ belong to the class $\W_3$.

According to \eqref{2.5}, we have
\begin{equation}\label{3.6}
    g\bigl([X_i,X_j],X_i\bigr)=g\bigl([X_i,X_j],X_j\bigr)=0.
\end{equation}
Then the following decomposition is valid
\begin{equation}\label{3.7}
\begin{array}{ll}
    [X_1,X_2]= C_{12}^3 X_3 +C_{12}^4 X_4,\qquad & [X_2,X_3]= C_{23}^1 X_1 +C_{23}^4
    X_4,\\[4pt]
    [X_1,X_3]= C_{13}^2 X_2 +C_{13}^4 X_4,\qquad & [X_2,X_4]= C_{24}^1 X_1
    +C_{24}^3 X_3,\\[4pt]
    [X_1,X_4]= C_{14}^2 X_2 +C_{14}^3 X_3,\qquad & [X_3,X_4]= C_{34}^1 X_1
    +C_{34}^2 X_2.\\[4pt]
\end{array}
\end{equation}

Now we apply again \eqref{2.5} using \eqref{3.7}. So we obtain
\begin{equation}\label{3.8}
\begin{array}{ll}
    [X_1,X_2]= \lambda_1 X_3 +\lambda_2 X_4,\qquad & [X_2,X_3]= \lambda_1 X_1
    +\lambda_3 X_4,\\[4pt]
    [X_1,X_3]= -\lambda_1 X_2 +\lambda_4 X_4,\qquad & [X_2,X_4]= \lambda_2 X_1
    -\lambda_3 X_3,\\[4pt]
    [X_1,X_4]= -\lambda_2 X_2 -\lambda_4 X_3,\qquad & [X_3,X_4]= \lambda_4 X_1
    +\lambda_3 X_2,\\[4pt]
\end{array}
\end{equation}
where $\lambda_1=C_{12}^3$, $\lambda_2=C_{12}^4$,
$\lambda_3=C_{23}^4$, $\lambda_4=C_{13}^4$. We verify immediately
that the Jacobi identity is satisfied in this case.

Let the conditions \eqref{3.8} be satisfied for a Riemannian
almost product manifold $(G,P,g)$ with structure $P$ and metric
$g$, determined by \eqref{3.1} and \eqref{2.2}, respectively. Then
we verify directly that $g$ is a Killing metric.

Therefore, the following theorem is valid.
\begin{thm}\label{thm-3.1}
    Let $(G,P,g)$ be a 4-dimensional Riemannian
almost product manifold, where $G$ is the connected Lie group with
an associated Lie algebra, determined by a global basis $\{X_i\}$
of left invariant vector fields, and $P$ and $g$ are the almost
product structure and the Riemannian metric, determined by
\eqref{3.1} and \eqref{2.2}, respectively. Then $(G,P,g)$ is a
$\W_3$-manifold with a Killing metric $g$ iff $G$ belongs to the
4-parametric family of Lie groups, determined by \eqref{3.8}.
\end{thm}

From this point on, until the end of this section we shall
consider the Riemannian almost product manifold $(G,P,g)$
determined by the conditions of \thmref{thm-3.1}.

Using \eqref{3.4}, \eqref{3.8}, \eqref{3.1} and \eqref{3.2}, we
obtain the following nonzero components of the tensor $F$:
\begin{equation}\label{3.9}
\begin{array}{l}
    F_{211}=-F_{233}=2F_{134}=2F_{323}=-2F_{112}=-2F_{314}=\lambda_1,\\[4pt]
    F_{144}=-F_{122}=2F_{212}=2F_{423}=-2F_{234}=-2F_{414}=\lambda_2,\\[4pt]
    F_{322}=-F_{344}=2F_{214}=2F_{434}=-2F_{223}=-2F_{412}=\lambda_3,\\[4pt]
    F_{433}=-F_{411}=2F_{141}=2F_{321}=-2F_{132}=-2F_{334}=\lambda_4.\\[4pt]
\end{array}
\end{equation}
The other nonzero components of $F$ are obtained from the property
$F_{ijk}=F_{ikj}$.

Let $F$ be the Nijenhuis tensor on $(G,P,g)$, i.e.
\[
    N_{ij}=[X_i,X_j]+P[PX_i,X_j]+P[X_i,PX_j]-[PX_i,PX_j].
\]
According to \eqref{3.1} and \eqref{3.8}, for the square norm
$\nN=N_{ik}N_{js}g^{ij}g^{ks}$ of $N$ we get
\begin{equation}\label{3.10}
    \nN=32\left(\lambda_1^2+\lambda_2^2+\lambda_3^2+\lambda_4^2\right).
\end{equation}

For the square norm of $\nabla P$, using \eqref{1.5}, \eqref{2.2}
and \eqref{3.3}, we obtain
\begin{equation}\label{3.11}
    \nP=4\left(\lambda_1^2+\lambda_2^2+\lambda_3^2+\lambda_4^2\right).
\end{equation}

From \eqref{2.7}, having in mind \eqref{2.2} and \eqref{3.8}, we
receive the following nonzero components of the curvature tensor
$R$:
\begin{equation}\label{3.12}
\begin{array}{ll}
    R_{1221}=\frac{1}{4}\left(\lambda_1^2+\lambda_2^2\right),\quad
    &
    R_{1331}=\frac{1}{4}\left(\lambda_1^2+\lambda_4^2\right),\\[4pt]
    R_{1441}=\frac{1}{4}\left(\lambda_2^2+\lambda_4^2\right),\quad
    &
    R_{2332}=\frac{1}{4}\left(\lambda_1^2+\lambda_3^2\right),\\[4pt]
    R_{2442}=\frac{1}{4}\left(\lambda_2^2+\lambda_3^2\right),\quad
    &
    R_{3443}=\frac{1}{4}\left(\lambda_3^2+\lambda_4^2\right),\\[4pt]
    R_{1341}=R_{2342}=\frac{1}{4}\lambda_1\lambda_2,\quad
    &
    R_{3123}=R_{4124}=\frac{1}{4}\lambda_3\lambda_4,\\[4pt]
    R_{1231}=R_{4234}=\frac{1}{4}\lambda_2\lambda_4,\quad
    &
    R_{2142}=R_{3143}=\frac{1}{4}\lambda_1\lambda_3,\\[4pt]
    R_{1241}=R_{3243}=-\frac{1}{4}\lambda_1\lambda_4,\quad
    &
    R_{2132}=R_{4134}=-\frac{1}{4}\lambda_2\lambda_3.\\[4pt]
\end{array}
\end{equation}
The other nonzero components of $R$ are obtained from the
properties $R_{ijks}=R_{ksij}$ and $R_{ijks}=-R_{jiks}=-R_{ijsk}$.

From \eqref{1.3}, having in mind \eqref{2.2}, we receive the
components $\rho_{ij}=\rho(X_i,X_j)$ of the Ricci tensor $\rho$.
The nonzero components of $\rho$ are:
\begin{equation}\label{3.13}
\begin{array}{c}
\begin{array}{ll}
    \rho_{11}=\frac{1}{2}\left(\lambda_1^2+\lambda_2^2+\lambda_4^2\right),\quad
    &
    \rho_{22}=\frac{1}{2}\left(\lambda_1^2+\lambda_2^2+\lambda_3^2\right),\\[4pt]
    \rho_{33}=\frac{1}{2}\left(\lambda_1^2+\lambda_3^2+\lambda_4^2\right),\quad
    &
    \rho_{44}=\frac{1}{2}\left(\lambda_2^2+\lambda_3^2+\lambda_4^2\right),\\[4pt]
\end{array}
\\[4pt]
\begin{array}{lll}
    \rho_{12}=\frac{1}{2}\lambda_3\lambda_4,\quad
    &
    \rho_{13}=-\frac{1}{2}\lambda_2\lambda_3,\quad
    &
    \rho_{14}=\frac{1}{2}\lambda_1\lambda_3,\\[4pt]
    \rho_{23}=\frac{1}{2}\lambda_2\lambda_4,\quad
    &
    \rho_{24}=-\frac{1}{2}\lambda_1\lambda_4,\quad
    &
    \rho_{34}=\frac{1}{2}\lambda_1\lambda_2.\\[4pt]
\end{array}
\end{array}
\end{equation}
The other nonzero components of $\rho$ are obtained from the
property $\rho_{ij}=\rho_{ji}$.

For the scalar curvature $\tau$, using \eqref{1.4}, we obtain
\begin{equation}\label{3.14}
    \tau=\frac{3}{2}\left(\lambda_1^2+\lambda_2^2+\lambda_3^2+\lambda_4^2\right).
\end{equation}

From \eqref{1.6}, having in mind \eqref{2.2}, \eqref{3.12},
\eqref{3.13} and \eqref{3.14}, we get for the Weyl tensor $W=0$.
Then $(G,P,g)$ is a conformally flat manifold.

For the sectional curvatures $k_{ij}=k(\alpha_{ij})$ of basic
2-planes $\alpha_{ij}=(X_i,X_j)$, according to \eqref{1.7},
\eqref{3.12} and \eqref{2.2}, we have:
\begin{equation}\label{3.15}
    \begin{array}{ll}
    k_{12}=\frac{1}{4}\left(\lambda_1^2+\lambda_2^2\right),\quad
    &
    k_{13}=\frac{1}{4}\left(\lambda_1^2+\lambda_4^2\right),\\[4pt]
    k_{14}=\frac{1}{4}\left(\lambda_2^2+\lambda_4^2\right),\quad
    &
    k_{23}=\frac{1}{4}\left(\lambda_1^2+\lambda_3^2\right),\\[4pt]
    k_{24}=\frac{1}{4}\left(\lambda_2^2+\lambda_3^2\right),\quad
    &
    k_{34}=\frac{1}{4}\left(\lambda_3^2+\lambda_4^2\right).\\[4pt]
\end{array}
\end{equation}

The obtained geometric characteristics of the considered manifold
are generalized in the following
\begin{thm}\label{thm-3.2}
    Let $(G,P,g)$ be the 4-dimensional Riemannian
almost product manifold where $G$ is the Lie group determined by
\eqref{3.8}, and the structure $P$ and the metric $g$ are
determined by \eqref{3.1} and \eqref{2.2}, respectively. Then
    \begin{enumerate}
    \renewcommand{\labelenumi}{(\roman{enumi})}
    \item
    $(G,P,g)$ is a locally symmetric $\W_3$-manifold with Killing metric $g$ and zero Weyl tensor;
    \item
    The nonzero components of the basic tensor $F$, the curvature tensor $R$ and the Ricci tensor $\rho$ are
    \eqref{3.9}, \eqref{3.12} and \eqref{3.13}, respectively;
    \item
    The square norms of the Nijenhuis tensor $N$ and $\nabla P$ are \eqref{3.10} and \eqref{3.11}, respectively;
    \item
    The scalar curvature $\tau$ and the sectional curvatures $k_{ij}$
    of the basic 2-planes are \eqref{3.14} and \eqref{3.15}, respectively.
    \end{enumerate}
\end{thm}

Let us remark that the 2-planes $\alpha_{13}$ and $\alpha_{24}$
are \emph{$P$-invariant 2-planes}, i.e.
$P\alpha_{13}=\alpha_{13}$, $P\alpha_{24}=\alpha_{24}$. The
2-planes $\alpha_{12}$, $\alpha_{14}$, $\alpha_{23}$,
$\alpha_{34}$ are \emph{totally real 2-planes}, i.e.
$\alpha_{12}\perp P\alpha_{12}$, $\alpha_{14}\perp P\alpha_{14}$,
$\alpha_{23}\perp P\alpha_{23}$, $\alpha_{34}\perp P\alpha_{34}$.
Then the equalities \eqref{3.15} imply the following
\begin{thm}\label{thm-3.3}
    Let $(G,P,g)$ be the 4-dimensional Riemannian
almost product manifold where $G$ is the Lie group determined by
\eqref{3.8}, and the structure $P$ and the metric $g$ are
determined by \eqref{3.1} and \eqref{2.2}, respectively. Then
    \begin{enumerate}
    \renewcommand{\labelenumi}{(\roman{enumi})}
    \item
    $(G,P,g)$ is of constant $P$-invariant sectional curvatures
    iff
    \[\lambda_1^2+\lambda_4^2=\lambda_2^2+\lambda_3^2;\]
    \item
    $(G,P,g)$ is of constant totally real sectional curvatures
    iff
    \[\lambda_1^2=\lambda_4^2,\qquad \lambda_2^2=\lambda_3^2\].
    \end{enumerate}
\end{thm}


\section{A Lie group as a pseudo-Riemannian almost product manifold
with Killing metric and nonintegrable structure}

In this section we consider a pseudo-Riemannian manifold $(G,P,g)$
with a metric $g$ determined by \eqref{2.3} and a structure $P$
defined as follows
\begin{equation}\label{4.1}
    PX_1=X_1,\quad PX_2=X_1,\quad PX_3=-X_3,\quad PX_4=-X_4.
\end{equation}
Obviously, $P^2=\id$. Moreover, \eqref{2.3} and \eqref{4.1} imply
\begin{equation}\label{4.2}
    g(PX_i,PX_j)=g(X_i,X_j).
\end{equation}
Therefore, $(G,P,g)$ is a pseudo-Riemannian almost product
manifold.

For the  manifold $(G,P,g)$ we propose that $g$ be a Killing
metric. Then $(G,P,g)$ is locally symmetric and the equalities
\eqref{2.6}, \eqref{2.7}, \eqref{3.3} and \eqref{3.4} are valid.

From \eqref{2.3}, \eqref{4.1} and \eqref{4.2} we obtain
\eqref{3.5}, i.e. $(G,P,g)$ is a $\W_3$-manifold.

Now, the equalities \eqref{3.6} and \eqref{3.7} are also
satisfied. According to \eqref{2.5}, from \eqref{3.7} we obtain
\begin{equation}\label{4.3}
\begin{array}{ll}
    [X_1,X_2]= \lambda_2 X_3 -\lambda_1 X_4,\qquad & [X_2,X_3]= -\lambda_2 X_1
    -\lambda_3 X_4,\\[4pt]
    [X_1,X_3]= \lambda_2 X_2 +\lambda_4 X_4,\qquad & [X_2,X_4]= \lambda_1 X_1
    +\lambda_3 X_3,\\[4pt]
    [X_1,X_4]= -\lambda_1 X_2 -\lambda_4 X_3,\qquad & [X_3,X_4]= -\lambda_4 X_1
    +\lambda_3 X_2,\\[4pt]
\end{array}
\end{equation}
where $\lambda_1=C_{24}^1$, $\lambda_2=C_{12}^3$,
$\lambda_3=C_{24}^3$, $\lambda_4=C_{13}^4$. We verify immediately
that the Jacobi identity is satisfied in this case.

Let the conditions \eqref{4.3} be satisfied for a
pseudo-Riemannian almost product manifold $(G,P,g)$ with structure
$P$ and metric $g$ determined by \eqref{4.1} and \eqref{2.3},
respectively. Then we verify directly that $g$ is a Killing
metric.

Therefore, the following theorem is valid.
\begin{thm}\label{thm-4.1}
    Let $(G,P,g)$ be a 4-dimensional pseudo-Riemannian
almost product manifold, where $G$ is the connected Lie group with
an associated Lie algebra, determined by a global basis $\{X_i\}$
of left invariant vector fields, and $P$ and $g$ are the almost
product structure and the pseudo-Riemannian metric, determined by
\eqref{4.1} and \eqref{2.3}, respectively. Then $(G,P,g)$ is a
$\W_3$-manifold with a Killing metric $g$ iff $G$ belongs to the
4-parametric family of Lie groups, determined by \eqref{4.3}.
\end{thm}

From this point on, until the end of this section we shall
consider the pseudo-Riemannian almost product manifold $(G,P,g)$
determined by the conditions of \thmref{thm-4.1}.

In an analogous way of the previous section, we get some geometric
characteristics of $(G,P,g)$.

We obtain the following nonzero components of the tensor $F$:
\begin{equation}\label{4.4}
\begin{array}{ll}
    F_{124}=-F_{214}=\lambda_1, \qquad &
    F_{213}=-F_{123}=\lambda_2,\\[4pt]
    F_{423}=-F_{324}=\lambda_3, \qquad &
    F_{314}=-F_{413}=\lambda_4.\\[4pt]
\end{array}
\end{equation}
The other nonzero components of $F$ are obtained from the
properties $F_{ijk}=F_{ikj}$.

The square norms of the Nijenhuis tensor $N$ and $\nabla P$ are
respectively:
\begin{equation}\label{4.5}
    \nN=24\left(\lambda_1^2+\lambda_2^2-\lambda_3^2-\lambda_4^2\right),
\end{equation}
\begin{equation}\label{4.6}
    \nP=-4\left(\lambda_1^2+\lambda_2^2-\lambda_3^2-\lambda_4^2\right).
\end{equation}

The nonzero components of the curvature tensor $R$ and the Ricci
tensor $\rho$ are respectively:
\begin{equation}\label{4.7}
\begin{array}{ll}
    R_{1221}=-\frac{1}{4}\left(\lambda_1^2+\lambda_2^2\right),\quad
    &
    R_{1331}=\frac{1}{4}\left(\lambda_2^2-\lambda_4^2\right),\\[4pt]
    R_{1441}=-\frac{1}{4}\left(\lambda_1^2-\lambda_4^2\right),\quad
    &
    R_{2332}=\frac{1}{4}\left(\lambda_2^2-\lambda_3^2\right),\\[4pt]
    R_{2442}=\frac{1}{4}\left(\lambda_1^2-\lambda_3^2\right),\quad
    &
    R_{3443}=\frac{1}{4}\left(\lambda_3^2+\lambda_4^2\right),\\[4pt]
    R_{1341}=R_{2342}=-\frac{1}{4}\lambda_1\lambda_2,\quad
    &
    R_{2132}=-R_{4134}=\frac{1}{4}\lambda_1\lambda_3,\\[4pt]
    R_{1231}=-R_{4234}=\frac{1}{4}\lambda_1\lambda_4,\quad
    &
    R_{2142}=-R_{3143}=\frac{1}{4}\lambda_2\lambda_3,\\[4pt]
    R_{1241}=-R_{3243}=\frac{1}{4}\lambda_2\lambda_4,\quad
    &
    R_{3123}=R_{4124}=\frac{1}{4}\lambda_3\lambda_4;\\[4pt]
\end{array}
\end{equation}
\begin{equation}\label{4.8}
\begin{array}{c}
\begin{array}{ll}
    \rho_{11}=-\frac{1}{2}\left(\lambda_1^2+\lambda_2^2-\lambda_4^2\right),\quad
    &
    \rho_{22}=-\frac{1}{2}\left(\lambda_1^2+\lambda_2^2-\lambda_3^2\right),\\[4pt]
    \rho_{33}=\frac{1}{2}\left(\lambda_2^2-\lambda_3^2-\lambda_4^2\right),\quad
    &
    \rho_{44}=\frac{1}{2}\left(\lambda_1^2+\lambda_3^2-\lambda_4^2\right),\\[4pt]
\end{array}
\\[4pt]
\begin{array}{lll}
    \rho_{12}=-\frac{1}{2}\lambda_3\lambda_4,\quad
    &
    \rho_{13}=\frac{1}{2}\lambda_1\lambda_3,\quad
    &
    \rho_{14}=\frac{1}{2}\lambda_2\lambda_3,\\[4pt]
    \rho_{23}=\frac{1}{2}\lambda_1\lambda_4,\quad
    &
    \rho_{24}=\frac{1}{2}\lambda_2\lambda_4,\quad
    &
    \rho_{34}=-\frac{1}{2}\lambda_1\lambda_2.\\[4pt]
\end{array}
\end{array}
\end{equation}
The other nonzero components of $R$ and $\rho$ are obtained from
the properties $R_{ijks}=R_{ksij}$, $R_{ijks}=-R_{jiks}=-R_{ijsk}$
and $\rho_{ij}=\rho_{ji}$.

The scalar curvature is
\begin{equation}\label{4.9}
    \tau=-\frac{3}{2}\left(\lambda_1^2+\lambda_2^2-\lambda_3^2-\lambda_4^2\right).
\end{equation}

We get for the Weyl tensor that $W=0$. Then $(G,P,g)$ is a
conformally flat manifold.

The sectional curvatures $k_{ij}=k(\alpha_{ij})$ of basic 2-planes
$\alpha_{ij}=(X_i,X_j)$ are:
\begin{equation}\label{4.10}
    \begin{array}{ll}
    k(\alpha_{13})=-\frac{1}{4}\left(\lambda_2^2-\lambda_4^2\right),\quad
    &
    k(\alpha_{24})=-\frac{1}{4}\left(\lambda_1^2-\lambda_3^2\right),\\[4pt]
    k(\alpha_{12})=-\frac{1}{4}\left(\lambda_1^2+\lambda_2^2\right),\quad
    &
    k(\alpha_{14})=-\frac{1}{4}\left(\lambda_1^2-\lambda_4^2\right),\\[4pt]
    k(\alpha_{23})=-\frac{1}{4}\left(\lambda_2^2-\lambda_3^2\right),\quad
    &
    k(\alpha_{34})=\frac{1}{4}\left(\lambda_3^2+\lambda_4^2\right).\\[4pt]
\end{array}
\end{equation}

Since $\alpha_{ij}=P\alpha_{ij}$ then all basic 2-planes are
$P$-invariant. It is used to check that now $(G,P,g)$ does not
accept constant $P$-invariant sectional curvatures.

The obtained geometric characteristics of the considered manifold
are generalized in the following
\begin{thm}\label{thm-4.2}
    Let $(G,P,g)$ be the 4-dimensional pseudo-Riemannian
almost product manifold where $G$ is the Lie group determined by
\eqref{4.3}, and the structure $P$ and the metric $g$ are
determined by \eqref{4.1} and \eqref{2.3}, respectively. Then
    \begin{enumerate}
    \renewcommand{\labelenumi}{(\roman{enumi})}
    \item
    $(G,P,g)$ is a locally symmetric conformally flat $\W_3$-manifold with Killing metric $g$;
    \item
    The nonzero components of the basic tensor $F$, the curvature tensor $R$ and the Ricci tensor $\rho$ are
    \eqref{4.4}, \eqref{4.7} and \eqref{4.8}, respectively;
    \item
    The square norms of the Nijenhuis tensor $N$ and $\nabla P$ are \eqref{4.5} and \eqref{4.6}, respectively;
    \item
    The scalar curvature $\tau$ and the sectional curvatures $k_{ij}$
    of the basic 2-planes are \eqref{4.9} and \eqref{4.10}, respectively.
    \end{enumerate}
\end{thm}

The last theorem implies immediately the following
\begin{cor}
 Let $(G,P,g)$ be the 4-dimensional pseudo-Riemannian
almost product manifold where $G$ is the Lie group determined by
\eqref{4.3}, and the structure $P$ and the metric $g$ are
determined by \eqref{4.1} and \eqref{2.3}, respectively. Then the
following propositions are equivalent:
   \begin{enumerate}
    \renewcommand{\labelenumi}{(\roman{enumi})}
    \item
    $(G,P,g)$ is an isotropic $P$-manifold;
    \item
    $(G,P,g)$ is a scalar flat manifold;
    \item
    The Nijenhuis tensor  is isotopic;
    \item
    The condition $\lambda_1^2+\lambda_2^2-\lambda_3^2-\lambda_4^2=0$ is
    valid.
    \end{enumerate}
\end{cor}


\bigskip

\textit{Dimitar Mekerov\\
University of Plovdiv\\
Faculty of Mathematics and Informatics
\\
Department of Geometry\\
236 Bulgaria Blvd.\\
Plovdiv 4003\\
Bulgaria
\\
e-mail: mircho@uni-plovdiv.bg}


\begin{thebibliography}{00}


\bibitem{Hel}
S.~Helgason, Differential geometry, Lie groups and symmetric
spaces, Academic Press, New York, 1978.



\bibitem{MaGrMe-4}
M.~Manev, K.~Gribachev, D.~Mekerov, On three-parametric Lie groups
as quasi-K\"ahler manifolds with Killing Norden metric. In: Topics
in Contemporary Differential Geometry, Complex Analysis and
Mathematical Physics, Eds. S. Dimiev and K. Sekigawa, World Sci.
Publ., Hackensack, NJ, 2007, 205-214.



\bibitem{StGr:connect}
M. Staikova, K. Gribachev, Cannonical connections and conformal
invariants on Riemannian almost product manifolds, Serdica Math.
J. {\bf 18} (1992), 150-161.


\end{thebibliography}
\end{document}